\documentclass[12pt]{article}
\usepackage{latexsym}
\usepackage{amsfonts}
\title{Genetic algorithms and the  Andrews-Curtis conjecture}
\author{Alexei D. Miasnikov}
\date{}
\newtheorem{theorem}{Theorem}

\newtheorem{lemma}{Lemma}

\newtheorem{proposition}{Proposition}

\makeatletter

\def\classification{\@ifnextchar [{\@xfootnotenext}%
   {\begingroup\let\protect\noexpand
      \xdef\@thefnmark{}\endgroup
    \@footnotetext}}

\begin{document}

\maketitle
\classification {{\it 1991 Mathematics
Subject Classification:} Primary 20E05, 20F05, 68T05; Secondary  57M05,
57M20. }

\begin{abstract}
\noindent  The Andrews-Curtis conjecture claims that every balanced
presentation of the trivial group can be
transformed into the trivial presentation
 by a finite
sequence of ``elementary  transformations" which are Nielsen
transformations together with an arbitrary conjugation of a relator.
It is believed that the Andrews-Curtis conjecture is false;
however, not so many possible counterexamples are known.
It is not a trivial matter to verify whether the conjecture holds
for a given balanced  presentation or not. The purpose of this paper is
to describe some non-deterministic  methods, called Genetic Algorithms,
designed to test the validity of the Andrews-Curtis conjecture.
Using such algorithm we have been able to prove  that all known (to us)
balanced presentations of the trivial group where the total length
of the relators is at most 12  satisfy  the conjecture.
In particular, the Andrews-Curtis
conjecture holds for the presentation \[<x,y|x y x = y x y, x^2 = y^3>\]
which was one of the well known  potential counterexamples.

\end{abstract}
\section{Introduction}

One of the difficulties in working with finitely presented groups is
the  fact that a lot of problems about them are unsolvable. Even if algorithms
exist, many of them are exponential or super exponential in nature making
it very unlikely for them to produce  results in an
acceptable period of time.
That is why the idea of applying non-deterministic methods to solve such
complex problems seems to be promising. In this paper we show that
Genetic Algorithms work fairly  well when applied to the
famous Andrews-Curtis conjecture.

The problem that we will be concerned with here, now termed
the \textit{Andrews-Curtis conjecture} was raised
by  J.J.Andrews and M.L.Curtis in their paper \cite{AC}. It
is of interest in topology as well as in group theory.  A good
group-theoretical survey was given by R.G.Burns and
Olga Macedonska \cite{BM}.\\

\textbf{Andrews-Curtis conjecture} \cite{AC}.
\textit{
If $$<x_{1},...,x_{n};r_{1},...,r_{n}>$$
 is a presentation of the trivial group,
can  it be converted to the trivial presentation
$$<x_{1},...,x_{n};x_1,\dots,x_n>$$
 by a finite sequence of the
following transformations: \\
(AC1) replace some $r_{i}$ by $r_{i}^{-1}$, \\
(AC2) replace some $r_{i}$ by $r_{i}r_{j}, j\neq i$, \\
(AC3) replace some $r_{i}$ by $wr_{i}w^{-1}$, where $w$ is any word in the
generators,\\
(AC4) re-order the relators; \\
(AC5) introduce a new generator $y$ and the relator $y$ or delete
a generator $y$ and relator $y$.
}
\\

We say that two presentations $G$ and $G'$ are Andrews-Curtis equivalent
(AC-equivalent) if one of them can be obtained from the other by
applying a finite sequence of transformations of the type $(AC1)-(AC4)$.

There is a weak form of the Andrews-Curtis conjecture when one is
allowed to use the transformations ($AC1$)-($AC5$). In this case we
say that the presentations $G$ and $G'$ are stably  AC-equivalent.
We will focus on the strong form of the conjecture because all of our results
 were obtained by using the transformations $(AC1)-(AC4)$.

There are series of presentations of the trivial group which are
known to satisfy the conjecture. Some of them are given in sections 4 and 5.
However, it is generally believed that the Andrews-Curtis conjecture is false.
There are a number of presentations which are considered as
potential counterexamples.

Here are some well known presentations of the trivial group that have
been put forward previously as possible  counterexamples:\\
$(1) <x,y,z;y^{-1} x y = x^2,z^{-1}yz=y^2,x^{-1}zx=z^2>;\\$
$ (2) <x,y;x^{-1}y^2x=y^3,y^{-1}x^2y=x^3>;\\$
$ (3) <x,y;x^2 = y^3,xyx=yxy>.$
 \\

The first two presentations have been known for almost 20 years;
for a discussion
see the paper by R.G.Burns and Olga Macedonska \cite{BM}.
The presentation (3) has been known for fifteen years.
It is the shortest one in the series $ <x,y;x^n = y^{n+1},xyx=yxy>,$
$ n \geq 2$ which appeared in \cite{AK} in 1985.

The purpose of this paper is to describe new genetic algorithms
designed to test the validity of the Andrews-Curtis conjecture and
to present results on the conjecture obtained with the help of
these algorithms. The main result (see section \ref{sec:results})
tells one that  all balanced presentations
of the trivial group that we are aware of where the total length
of defining relators is at most 12 satisfy the conjecture.

In particular the presentation (3)  in fact  satisfies the conjecture.
We found a chain of 21 transformations that reduces the presentation
(3) to the trivial presentation ( see section \ref{sec:results}).
One of the most general series of balanced presentations of the trivial
group is contained in the recent paper of C.F.Miller and P.E.Schupp \cite{MS}:
$<x,y;x^{-1}y^nx=y^{n+1},x=w>$ ,
where $n \geq 1$, and $w$ is a word in $x$ and $y$ with exponent
sum 0 in $x$.

In section 4 we describe some experiments which show the surprisingly high
capabilities of genetic algorithms. Indeed we could not find an example
of a balanced presentation of the trivial group which is known to be
AC-equivalent to the trivial presentation where our algorithm
couldn't find a solution. Of course we considered only examples of a
reasonable size.
Some basic terminology and the algorithm itself will be introduced
in sections 2 and 3.

In a little more detail, then, we want to devise algorithms which carry
out a sequence of  transformations which
transform a  given presentation of the trivial group into the trivial one.
From the point of view of computation, neither
total enumeration, nor random search can be effectively applied here.
In order to enumerate all sequences of transformations we have to
implement some adjustments to the original procedure.
We  substitute the transformation $(AC3)$ by two new transformations: \\
($AC3'$) replace some $r_{i}$ by $x_{j}r_{i}x_{j}^{-1}$,\\
($AC3''$) replace some $r_{i}$ by $x_{j}^{-1}r_{i}x_{j}$.\\
It is clear that the transformation $(3)$ can be obtained by a sequence of
transformations  ($AC3'$) and ($AC3''$).  We  mention
here that all of our solutions were obtained without using the
transformation ($AC5$) and we were unable to improve efficiency  by
adding this operation. Transformation ($AC4$)
is not important here and we will ignore it for the moment.

Now, disregarding the transformations ($AC4$) and ($AC5$),
 we have a total of $3n^2$ elementary
transformations, where $n$ is the number of generators (relators) in
the presentation. In order to check all sequences of transformations of length $k$, we need
to produce $(3n^2)^{k}$ chains of transformations. Even with $n=2$ it will
be hard to get to a fairly long sequence because we need to check all
sequences starting from the empty one.
For example to arrive at the solution mentioned above,
we would have to enumerate  $12^{21} \approx 4.6*10^{22} $ sequences of transformations, which is
impossible even for a powerful computer.

The number of computations can be reduced at expense of memory by
using the following approach. Suppose we want do check weather
a given presentation, say $G$, can be reduced to the
trivial presentation of the trivial group by a sequence of
transformations ($AC1$)-($AC3$) and the length of such sequence is not greater than  21.
First, we produce and store in memory all presentations obtained from
the trivial one by transforming it with the sequences of transformations with the
length at most 10. The next step is to apply chains of
transformations with the length less or equal 11 to the
presentation $G$. If at least one of the presentations obtained from $G$
lies in the ball saved during the previous step then
the presentation $G$ satisfies the Andrews-Curtis conjecture.

This approach reduces the number of computations to $\approx 12^{11}$.
Unfortunately it requires  a lot of a fast access memory (more than 60 Gb in
our example) in order to be effective. 

We ran the total enumeration
procedure on a machine with an Alpha processor running at 500 mhz
and only managed to test sequences of length up to seven in two days.
Random enumerations are also unlikely to give any
results because of the very small probability of finding the right sequence.
\section{Introduction to Genetic Algorithms}
\textit{Genetic algorithms} provide a method of searching for the optimal
solution in a space of all possible solutions corresponding to a
specific question. Terminology as well as the idea itself come from
biology. The first Genetic Algorithm goes back to J.H.Holland \cite{Hol}
in 1975.
This theory has been further developed and now Genetic Algorithms are widely
used to solve search and optimization problems. Genetic Algorithms may
vary in their implementation. There are many different approaches that are used
to solve particular problems. More detailed descriptions of Genetic
Algorithms and their applications can be found in the book by M.Mitchell
\cite{Mit}.
Here is a brief description of the major
components of the standard Genetic Algorithm.

Genetic Algorithms handle some selection of possible solutions -
\textit{population of solutions}. Each solution is encoded in a
specific representation called a
\textit{chromosome}. Originally, solutions were encoded into strings of bits.
Although binary codage is seen as a standard representation, usually
it is more convenient and effective to use problem-specific representations.
The algorithm starts by choosing an initial population. Generally, the first
population is generated randomly to provide a higher diversity of possible
solutions. However, some problems may require different initialization.

Members of a population are compared by using a \textit{fitness function}.
Each member (chromosome) has an associated value corresponding to the
fitness of the possible solution it represents. This value indicates how close
this chromosome is to an optimal solution.
 Usually a \textit{threshold value} or \textit{termination value} is
chosen as an
upper bound for the fitness values. The algorithm terminates when the fitness value
of one of the solutions reaches the threshold value. This solution is usually
optimal.

Populations evolve by applying \textit{reproduction operators} on specifically
chosen chromosomes. There are two basic forms of reproduction:
\textit{recombination} (or \textit{crossover}) and \textit{mutation}.
 Recombination is
the process of producing a child from two parents. The most
common reproduction operator is N-point crossover. Usually N-point
crossover takes two chromosomes and randomly cuts them into N+1 segments.
Children are produced by recombining alternatively chosen segments from their
parents. Mutation is a simple operator which randomly changes population
members. It is often seen as a background operator to maintain
diversity in the population and avoid premature convergence.

\textit{Selection} is a method that chooses members of the population
according to their fitness functions. Usually random selection methods
are used. The closer its fitness value is to the optimal one,
the more chance a
member has to be selected. Selection is responsible for selecting
chromosomes for  reproduction.
\begin{figure}[t]
\begin{verbatim}
BEGIN  /* genetic algorithm*/
  Generate initial population;
  Compute fitness values for all members of the population;

  WHILE NOT fitness values have reached threshold DO

     Select members from the current population for the
     recombination;

     Create new members by applying recombination and/or
     mutation to the selected members;

     Generate a new population by replacing members of the
     current population by the new ones according to the
     defined Replacement method;

  END WHILE
END
\end{verbatim}
\caption{Structure of the Genetic Algorithm.}
\end{figure}

The last  component of the Genetic Algorithm is \textit{Replacement}.
It is responsible for determining which of the current members of the
population, if any, should be replaced by the new ones. There are many
different approaches to produce replacement. Sometimes all members of
the current population are replaced by offsprings. Sometimes, in order
to increase convergence, only the best members (with highest fitness values)
among offsprings and parents are chosen to produce the new population.
The latter is called \textit{Strong Elitist Selection}.

Once all components of the Genetic Algorithm have been defined, it evolves
according to the basic structure which is  shown in Figure 1.

\section{Description of the algorithm}
First we have to define the space of all possible solutions for our
genetic algorithm to search in. We search for a sequence
of Andrews-Curtis transformations which transform our presentation
\[G = <x_{1},...,x_{n};r_{1},...,r_{n}>\] into another presentation
$G' = <x_{1},...,x_{n};r'_{1},...,r'_{n}>$ which satisfies some chosen
criteria encoded in the fitness function.

\subsection{Members of the population}
There are four elementary transformations $(AC1)-(AC4)$ that
are used in the algorithm. Since the transformation $(AC5)$ is
not involved, the set of generators $\{x_1,...,x_n\}$ does not
change during the transformation process. It follows, that in fact
we are transforming only the tuple of relators $(r_1,...,r_n)$. To
this end we say that tuples $U$ and $V$ from $F^n$ are AC-equivalent
(denoted $U\sim_{AC}V$ )
if $V$ can be obtained from $U$ by a finite sequence of transformations
$(AC1)-(AC4)$, where $F$ denotes the free group on $x_1,...,x_n$ .

We  define a set of elementary transformations $T$ which
 consists of the following:
\begin{enumerate}
\item Andrews-Curtis transformations $(AC1)-(AC4)$.

\item Whitehead automorphisms : \\
$
x_{i} \rightarrow x_{i}^{-1},\: x_{l} \rightarrow x_{l}, \\
x_{i} \rightarrow x_{j}^{\pm 1}x_{i},\: x_{l} \rightarrow x_{l}, \\
x_{i} \rightarrow x_{i}x_{j}^{\pm 1},\: x_{l} \rightarrow x_{l},\\
x_{i} \rightarrow x_{j}^{-1}x_{i}x_{j},\: x_{l} \rightarrow x_{l},\\
$
 where $i \neq j$ and  $i \neq l$.

\item Auxiliary transformations which are a special case of $(AC3)$: \\
replace some $r_{i}$ by $x_{j}r_{i}x_{j}^{-1}$,\\
replace some $r_{i}$ by $x_{j}^{-1}r_{i}x_{j}$,\\
produce a random cyclic permutation of one of the $r_{i}$.
\end{enumerate}

For any $n$-tuple $U=(r_1,...,r_n)$ and an elementary transformation
$t$ denote by $Ut$ the result of the transformation of $U$ by $t$.
If $ \varphi$ is one of the Whitehead automorphisms then we define
the effect of $\varphi$ on $U$ to be $(r_1^{\varphi},...,r_n^{\varphi})$
which we denote by $U^{\varphi}$.

\begin{lemma}
Let $U=(r_1,...,r_n)$ be an $n$-tuple of elements from $F=<x_1,...,x_n>$.
If the tuple $U$ can be transformed to the tuple of
 generators $X=(x_1,...,x_n)$  by a
sequence of transformations $T$ then $U$ is AC-equivalent to $X$.

\end{lemma}
\textit{Proof.} If $t$ is any transformation of type $(AC1),\;(AC2)$ or $(AC4)$
and $\varphi$ is a Whitehead automorphism, then
\begin{equation}
\label{E1}
U^{\varphi}t = (Ut)^{\varphi};
\end{equation}
if $t$ is a transformation of  type $(AC3)$, say
$(r_1,...,r_n)t = (r_1,...,r_i^{w},...,r_n)$ then
\begin{equation}
\label{E2}
U^{\varphi}t = (Us)^{\varphi},
\end{equation}
where the transformation $s$ conjugates $r_i$ by $w^{\varphi^{-1}}$.

Suppose now that the n-tuple $U$ can be transformed into $(x_1,...,x_n)$
by a finite sequence of transformations $T$. Then it follows from (\ref{E1})
and (\ref{E2}) that  we can find transformations $t_1,...,t_m$
of type $(AC1)-(AC4)$ and  Whitehead automorphisms
$\varphi_1,...,\varphi_k$ such that
\[(Ut_1...t_m)^{\varphi_1...\varphi_k} = (x_1,...,x_n).\]
It follows that
\[(Ut_1...t_m)^{\varphi_1...\varphi_{k-1}} =
(x_1,...,x_n)^{\varphi_k^{-1}}.\]

Observe that $(x_1,...,x_n)^{\varphi_k^{-1}} = (x_1,...,x_n)s$, where
$s$ is a transformation of type $(AC1)-(AC3)$. Therefore
\[(Ut_1...t_m)^{\varphi_1...\varphi_{k-1}}s^{-1} = (x_1,...,x_n).\]

It follows then from (\ref{E1}) and (\ref{E2}) that
\[(Ut_1...t_m)^{\varphi_1...\varphi_{k-1}}s^{-1} =
(Ut_1...t_m t_{m+1})^{\varphi_1...\varphi_{k-1}}\]
where $t_{m+1}$ is a transformation of type $(AC1)-(AC3)$.
In other words we find that
\[(Ut_1...t_m t_{m+1})^{\varphi_1...\varphi_{k-1}} = (x_1,...,x_n).\]
Iterating this procedure we can then deduce that
\[Ut_1...t_{m+k} = (x_1,...,x_n)\]
with $t_1,...,t_{m+k}$ transformations of type $(AC1)-(AC3)$. \\
This proves Lemma 1.

Since the transformation $(AC4)$ can be applied
at the very end, we avoid its use here.

The next result is that we have, in all,  a set $T'$, say,
 of  $8n^2 - 2n$ transformations to work with, i.e. all of the
 transformations of type $T$ except for  $(AC4)$.

The members of our population are represented by finite sequences of
transformations $t_i \in T'$:
\[member\;of\;the\;population\;=\;t_1,...,t_l,\;(t_i \in T').\]

Each such member represents a possible sequence which
transforms an n-tuple $(r_1,...,r_n)$ of relators into $(x_1,...,x_n)$.

\subsection{Fitness functions and termination conditions}
Implementation of the fitness function depends on the criteria by which we
will evaluate members of the population. Several different approaches were
implemented. Each of them is designed for slightly different
purposes and allows us to answer different questions.
Now, let $m=t_1,...,t_k$ be a member of our population. Suppose
the sequence $m$ transforms the given presentation
\[G=<x_1,...,x_n;r_1,...,r_n>\] into a presentation
\[G' = <x_{1},...,x_{n};r'_{1}, ...,r'_{n}>.\]
Below we discuss different criteria (fitness functions) which we used
for evaluating members of the population.

\subsubsection{Leading coordinates}
In order to prove that a presentation satisfies
the Andrews-Curtis conjecture
it is enough to prove that this presentation can be transformed by the
transformations $T'$ into
a presentation in which  $n-1$ relators form a part of a basis of $F_n$.
 This applies, in particular in the case $n=2$
which means that all we need to do is to transform one of the
relators into a primitive
element. Priority should be given to the transformations that reduce lengths
of relators. Also we don't care about the longest of them as we need to convert
only $n-1$ relators into primitive elements.

To check whether a new relator $r'_{i}$ is a primitive element we can apply
Whitehead's algorithm to reduce the length of the relator. If a relator
can be transformed to a generator then this relator is a primitive element.

A fitness function can be implemented as a function which depends on the
lengths of the relators. We denote the length of a word $w$ by $|w|$.

Let $R'=\{r'_{1},...,r'_{n}\}$. Choose one of the relators from $R'$ of
 maximal length. Denote all others by $s_1,...,s_{n-1}$ and put
$R'_{min}=\{s_1,...,s_{n-1}\}$.

Then the fitness value of a
member of a given population that transforms a given presentation $G$ into
its AC-equivalent presentation $G'$ can be defined as:
\[Fit1=\sum_{i=1}^{n-1}|s_{i}|,\]
where $s_{i} \in R'_{min} $.
Although it is customary in the application of genetic algorithms to define fitness
functions so that the optimal values are always  maximum values,
in our implementation the fittest member of the population  has the minimal value of $Fit1$.

The algorithm may terminate after a sequence of transformations from $T'$
transforms $\{r_{1},...,r_{n}\} \rightarrow \{ r'_{1},...,r'_{n}\} $ such that
$R'_{min}$  consists of only primitive elements. In order to see whether this
is the case, we apply Whitehead's algorithm to transform
 $\{s_{1},...,s_{n-1}\}$
into a set of relators $\{s'_{1},...,s'_{n-1}\}$:
\[\{s_{1},...,s_{n-1}\} \rightarrow \{s'_{1},...,s'_{n-1}\}, \]
the sum of whose lengths is minimal.
Now the algorithm will terminate when the following termination condition
holds:
\[\sum_{i=1}^{n-1}|s'_{i}|=n-1.\]

\subsubsection{Total length of relators}
There is another variation of a fitness function that appears  to be
very useful. Generally it is used when the algorithm is required
to produce the whole sequence of transformations $T$ that transforms the
relators $\{r_{1},...,r_{n}\} $ into the generators $\{x_{1},...,x_{n}\} $.
The idea is the same, only instead of reducing lengths of $n-1$ relators, we
 try to reduce the lengths of
all of them. The fittest transformation, again, should
have the minimal fitness, where the fitness function is given
by the following formula:
\[Fit2 = \sum_{i=1}^{n} |r'_{i}|,\]
where $r'_{i} \in R'$.

The termination value is taken to be $n$,
 which means that the algorithm will
terminate if each of the $r'_{i}$ has length 1.

\subsubsection{The shortest solution}
When it is known that a presentation satisfies the Andrews-Curtis conjecture,
sometimes it is  useful to produce not an arbitrary sequence of transformations
from the given presentation to the trivial presentation,
but instead one involving
the fewest number of such transformations. With this in mind,
let $t=t_{1},...,t_{k}$ be a sequence of transformations from $T'$; we term
$k$ the length of the sequence $t$ and define then $|t|=k$.
We  apply a penalty trick to the previously defined $Fit2$:
\[Fit3 = Fit2 + \frac{k}{m},\]
where  $m$ is a positive integer parameter defined by the user. It gives
priority to transformations with length less then $m$.
Longer chains will
have a smaller chance to be selected for recombination or to stay
 in the population.

The termination value must be adjusted, too. When checking for
the termination condition we have to keep in mind the penalty $\frac{k}{m}$.

Algorithms with fitness function $Fit3$ consume more time but produce
shorter sequences. On some examples it gave considerably shorter
answers.

\subsection{Reproduction operators}
\subsubsection{Crossover}
One point crossover was used in our genetic algorithms. It was applied with
some probability passed as a parameter to these algorithms.

To produce  recombination  two population members
$m_{1} = t_1,...,t_k$ and $m_{2}= s_1,...,s_l$
must be chosen in accordance with the selection method. Two random numbers
$p,\: 0 < p < k$ and $q,\: 0 < q < l$ define
the crossover points:
\[m_{1}=t_1,...,t_{p-1},t_p,...,t_k,\]
\[m_{2}=s_1,...,s_{q-1},s_q,...,s_l.\]

Two offsprings are obtained by recombination of  their parents' parts as
follows:
\[o_{1}=t_1,...,t_{p-1},s_q,...,s_l,\]
\[o_{2}=s_1,...,s_{q-1},t_p,...,t_k.\]

We have found that using two point and uniform crossovers did not produce
better results than the one-point crossover.

\subsubsection{Mutation}
Like crossover, we invoke mutation with some probability defined
as a parameter. Four elementary operations were introduced as mutations.
In order to explain, suppose that we want to apply a mutation to a
sequence of transformations
$t=t_{1},...,t_{k}$, where $t_{i} \in T'$. We then make use of the following
elementary operations:\\
\\
(M1) \textit{attach} a random transformation $s\in T'$ at the end of a
sequence of transformations  \[t\rightarrow t_{1},...,t_{k},s;\]
(M2)  \textit{insert} a random transformation $s\in T'$ into a
randomly chosen position $i$ in a sequence of transformations
\[t\rightarrow t_{1},...,t_{i-1},s,t_{i},..,t_{k};\]
(M3)  \textit{delete} a transformation in  a
randomly chosen position $i$ in a sequence of transformations
\[t\rightarrow t_{1},...,t_{i-1},t_{i+1},..,t_{k};\]
(M4) \textit{change} a transformation in a randomly chosen position to
a different one
\[t\rightarrow t_{1},...,t_{i-1},s,t_{i+1},..,t_{k}.\]

When applying mutation to the string $t$, only one of the operations
 $(M1)-(M4)$
is used.  Each of them has assigned a ``chance''
 value which defines the probability of an operation being applied.
This allows us to distribute the portions of applications of the mutation
operations.

Operator $(M1)$ is a special case of the more general operation
$(M2)$. It was specially introduced because we want to
 improve the fittest transformation while keeping its fitness unchanged.

\subsection{Selection and Replacement}
As a selection method the proportionate-based selection was chosen.
This is a probabilistic method of selection, and the probability $Pr(m)$
of the member $m$ to be selected is given by
\[ Pr(m)= \frac{F(m)}{\sum_{i=1}^p F(m_i)},\]
where $F(x)$ is a fitness value of the member $x$ scaled in the such
way that the higher value corresponds to the member which is
closest to the optimal.

The most classic fitness-proportionate selection ``Roulette
wheel selection'' method was implemented. Briefly, it simply assigns
to each possible solution a sector whose size is proportional to the
appropriate fitness measure. Then a random position on the wheel is
chosen. The selected chromosome is one which belongs to the sector containing
the chosen position of the wheel.

In order to increase diversity of
the fitness function evaluations the fitness measure is taken to be
the square of the fitness value.

All members of the population except the fittest one are replaced
by members of the new population. This type of replacement method produces
a more diverse population, but slows down the convergence of the procedure.
However replacement methods like ``Elitist selection''( members of
the parent population are replaced by  new ones only if they have
less fit values), which do converge more rapidly, do not seem
to yield better results. This is because of the "premature
convergence" - an effect when the algorithm falls into a
local minimum which is far from the optimal solution.

\section{Testing of capabilities}
It is believed that different parameters like crossover and mutation rates,
types of selection, replacement and recombination operators are
essential for better performance of Genetic Algorithms.
A lot of time was spent and many experiments were carried out so
as to figure out which
conditions provide the best results. It appears unlikely that there is a
common thread running through all of these algorithms.
Each of them  is very sensitive to its
parameters. In most cases the choice is problem-specific. 

We ran many experiments on the algorithms trying to find
parameter values that give the best performance.
This is probably not the
best way to proceed. I like the idea of algorithms which are
self-adjustable at run time.
This makes the algorithm more flexible and easy to use. Unfortunately
we have been unable to come up with effective methods for doing so.

Experiments showed that the algorithms work more effectively with 50 members
in the population. The probabilities of mutation and crossover were 
95\% and 85\% correspondingly and the mutations of type M1 had the greater chance
to appear.
 
Algorithms were tested on a list of examples of presentations known to
satisfy the Andrews-Curtis conjecture.
These include the following presentations: \\
(1)\(<x,y;x^{k} y^{l},x^{m} y^{n}>\), where $k,l,m,n$ - integers, such that
$kn-lm=\pm1$; \\
(2)\(<x,y;xc^{k}, yc^{l}>\), where $c \in [F_{2},F_{2}]$ and $k,l$ - any integers; \\
(3)\(<x,y;x^{-2}y^{-1}xy,bc>\), where $\{x,b\}$ is a free
basis for $F_{2}$ and $c \in [F_{2},F_{2}]$; \\
(4)$<x,y;r^s=r^2,s^r=s^2>$, where $a,b,r,s$ are such that $<a,b;r,s>$ is a
balanced presentation of the trivial group which satisfies the conjecture. \\

The presentations (1)-(3) were described in the survey of R.G.Burns and
O.Macedonska \cite{BM}. The last series (4) will be
discussed later in this paper (see the example after Proposition 1).

Examples of each kind of presentations (1)-(4) were included in the
testing list.  Here are some presentations that
we tested with our genetic algorithms: \\
1: $<a,b;(a^2 b^3)^{(a^3 b^4)} = (a^2 b^3)^2,(a^3 b^4)^{(a^2 b^3)} = (a^3 b^4)^2>$,
$ r = a^2 b^3,\:s = a^3 b^4$; \\
2: $<a,b;(a^3 b^4)^{(a^4 b^5)} = (a^3 b^4)^2,(a^4 b^5)^{(a^3 b^4)} = (a^4 b^5)^2>$,
$ r = a^3 b^4,\: s =  a^4 b^5$; \\
3: $<a,b;(a [a,b])^{(b [b,a])} = (a [a,b])^2,(b [b,a])^{(a [a,b])} = (b [b,a])^2>$,
$r = a [a,b],\: s = b [b,a]$; \\
4: $<a,b;(a^{-2} b^{-1} a b)^{(a b[a,b])}=(a^{-2} b^{-1} a b)^2,(a b[a,b])^{(a^{-2} b^{-1} a b)} =(a b[a,b])^2>$, $r = a^{-2} b^{-1} a b, \: s =  a b[a,b]$.

In each case the Andrews-Curtis conjecture was verified.

In addition,  we  tested the algorithms using
 automatically generated presentations. These we obtained  by
transforming the trivial presentation $<x,y;x,y>$  by
 a random sequence of transformations ($AC1$)-($AC3$). The results of this
experiment are shown in Table \ref{tab:table1}.

\begin{table}
\begin{tabular}{|l|l|l|l|}
\hline
Length of         & Average sum   & Average number of      &Number of  \\
a sequence of     & of relators'  & generations required   &examples   \\
transformations   & lengths       & to reach the result    &tested     \\
\hline
 10                &  13           & 20                     & 50 \\
\hline
 20                &  32           & 300                    & 50 \\
\hline
 30                &  78           & 2208                   & 50 \\
\hline
\end{tabular}
\caption{Experiments with automatically generated presentations}
\label{tab:table1}
\end{table}

To increase the diversity of the possible solutions we ran some
experiments with automatically generated presentations.
This experiments are very similar to those that we used for testing with
the only difference that sequences of random transformations $T$ were applied
to some presentation $G_t$ - a potential counterexample
(more often it was the presentation:  $<a,b;a b a = b a b, a^3 = b^4>$ ).
The new presentation is AC-equivalent to the original but has longer
relations. It is interesting to mention that in all
cases the algorithm converged to the original presentation $G_t$.

It is important to emphasize that in all tested cases our genetic algorithms
were able to find a sequence of transformations which reduced the tested
presentation into the trivial one. In fact we do not know of
any example of a presentation which is AC-equivalent to the
trivial presentation where our genetic algorithms failed to work.

\section{Results}
\label{sec:results}
Our major objective was to apply our genetic algorithms
 to some presentations that are considered to be likely
counter-examples to the Andrews-Curtis conjecture.

As we have mentioned already all interesting known balanced presentations
of the trivial group are given by two particular presentations:\\
$(1) <x,y,z;y^{-1} x y = x^2,z^{-1}yz=y^2,x^{-1}zx=z^2>;\\$
$(2) <x,y;x^{-1}y^2x=y^3,y^{-1}x^2y=x^3>,\\$
and two infinite series \\
$(3) <x,y;x^n = y^{n+1},xyx=yxy>,\; n \geq 2;$ \\
$(4) <x,y;x^{-1}y^nx=y^{n+1},x=w>$, where $n \geq 1$,
and $w$ is a word in $x$ and $y$ with exponent sum 0 in $x$.

The main result of this paper is the following theorem:
\begin{theorem}
All presentations from the series (3) and (4) with the total length
of relators at most 12 satisfy the Andrews-Curtis conjecture.
\end{theorem}
Altogether there are $273$ presentations with the total length $\leq$ 12 in
these series. Most of them are easily seen
to satisfy the conjecture, however some of them are hard to ``crack''.
Here we focus only on the presentations that are the hardest
to prove to satisfy the Andrews-Curtis conjecture.

We start with the well-known presentation $<a,b;a b a = b a b, a^2 = b^3>$
which is the only one in the series (3) with the total length of
relators $\leq 12$. We were able, using a genetic algorithm, to produce a sequence
of transformations from $T$ which transforms this presentation into
the trivial one.

Let \[G_1=<a,b; a b a = b a b, a^2 = b^3>\]
 and \[r_{0} = a^2 b^{-3},r_{1} = a b ab^{-1}a^{-1}b^{-1}.\]
We have the following chain of transformations: \\
$
1: r_{0} \rightarrow r_{0}^{-1} \Rightarrow r_{0} \rightarrow b^3 a^{-2}; \\
2: r_{1} \rightarrow r_{1} r_0 \Rightarrow r_{1} \rightarrow
a b a b^{-1} a^{-1} b^2 a^{-2}; \\
3: r_1 \rightarrow a^{-1} r_1 a \Rightarrow r_1 \rightarrow
b a b^{-1} a^{-1} b^2 a^{-1}; \\
4: r_{1} \rightarrow b^{-1} r_{1} b \Rightarrow r_{1} \rightarrow
a b^{-1} a^{-1} b^2 a^{-1} b;\\
5: r_{0} \rightarrow a^{-1} r_{0} a \Rightarrow r_{0} \rightarrow
a^{-1} b^3 a^{-1};\\
6: r_{0} \rightarrow r_{0} r_{1} \Rightarrow r_{0} \rightarrow
a^{-1} b^2 a^{-1} b^2 a^{-1} b; \\
7: r_{1} \rightarrow b^{-1} r_1 b \Rightarrow r_1 \rightarrow
b^{-1} a b^{-1} a^{-1} b^2 a^{-1} b^2; \\
8: r_{1} \rightarrow r_1 r_0 \Rightarrow r_1 \rightarrow
b^{-1} a b^{-1} a^{-1} b^2 a^{-1} b^2 a^{-1} b^2 a^{-1} b^2 a^{-1} b; \\
9: r_{1} \rightarrow b a^{-1} b r_1 b^{-1} a b^{-1} \Rightarrow r_1 \rightarrow
a^{-1} b^2 a^{-1} b^2 a^{-1} b^2 a^{-1} b;\\
10: r_{0} \rightarrow r_{0}^{-1} \Rightarrow r_{0} \rightarrow
b^{-1} a b^{-2} a b^{-2} a; \\
11: r_1 \rightarrow r_1 r_0 \Rightarrow r_1 \rightarrow a^{-1} b^2 ; \\
12: r_{0} \rightarrow r_0 r_1 \Rightarrow r_{0} \rightarrow
b^{-1} a b^{-2} a; \\
13: r_{0} \rightarrow r_0 r_1 \Rightarrow r_{0} \rightarrow
b^{-1} a; \\
14: r_{0} \rightarrow r_0 r_1 \Rightarrow r_{0} \rightarrow b; \\
15: r_{0} \rightarrow r_0^{-1} \Rightarrow r_{0} \rightarrow b^{-1}; \\
16: r_1 \rightarrow r_1 r_0 \Rightarrow r_1 \rightarrow a^{-1} b; \\
17: r_1 \rightarrow r_1 r_0 \Rightarrow r_1 \rightarrow a^{-1}; \\
18: r_1 \rightarrow r_1^{-1} \Rightarrow r_1 \rightarrow a; \\
19: r_0 \rightarrow r_0^{-1} \Rightarrow r_0 \rightarrow b.
$

The  transformation 9 requires the elementary transformations of type
$(AC3')$ and $(AC3'')$ adding two more transformations to the solution
 and increasing its length to 21.

Using  genetic algorithms we  were able to show that all the
presentations from the series (4) where the total length of relators
is at most 12 satisfy the conjecture as well. In order to do that we
enumerated all such presentations and checked
them with the algorithms.

It turns out that the presentations
\begin{equation}
\label{hard4}
<x,y;x^{-1} y^2 x = y^3, x = y^{\pm 1} x y^{\pm 1} x^{-1}>
\end{equation}
are the most interesting and hard to crack in the series (4).

Let
\[G_2=<x,y; y x y^{-1} x^{-2}, x^{-1} y^2 x y^{-3}>\]
 and denote
\[r_{0} = y x y^{-1} x^{-2},r_{1} = x^{-1} y^2 x y^{-3}.\] Then the
following chain of transformations reduces the presentation $G_2$
into the trivial group:
\\
$
1: r_1 \rightarrow x r_1 x^{-1} \Rightarrow
  r_1 \rightarrow y^2 x y^{-3} x^{-1};\\
2: r_1 \rightarrow  r_1^{-1} \Rightarrow
  r_1 \rightarrow x y^3 x^{-1} y^{-2};\\
3: r_1 \rightarrow y^{-1} r_1 y \Rightarrow
  r_1 \rightarrow y^{-1} x y^3 x^{-1} y^{-1};\\
4: r_1 \rightarrow  r_1 r_0 \Rightarrow
  r_1 \rightarrow y^{-1} x y^2 x^{-2};\\
5: x \rightarrow x^{-1} \Rightarrow
  r_0 \rightarrow y x^{-1} y^{-1} x^2,
  r_1 \rightarrow y^{-1} x^{-1} y^2 x^2;\\
6: r_1 \rightarrow  y r_1 y^{-1} \Rightarrow
  r_1 \rightarrow x^{-1}  y^2 x^2 y^{-1};\\
7: r_1 \rightarrow  r_1 r_0 \Rightarrow
  r_1 \rightarrow x^{-1}  y^2 x y^{-1} x^2;\\
8: r_1 \rightarrow x^2 r_1 x^{-2} \Rightarrow
  r_1 \rightarrow x y^2 x y^{-1};\\
9: r_1 \rightarrow  r_1 r_0 \Rightarrow
  r_1 \rightarrow x y x^2;\\
10: r_1 \rightarrow  r_1^{-1} \Rightarrow
  r_1 \rightarrow x^{-2} y^{-1} x^{-1};\\
11: y \rightarrow  y x^{-2} \Rightarrow
  r_0 \rightarrow y x^{-1} y^{-1} x^2,
  r_1 \rightarrow y^{-1} x^{-1};\\
12: r_1 \rightarrow  r_1^{-1} \Rightarrow
  r_1 \rightarrow x y;\\
13: r_0 \rightarrow  r_0^{-1} \Rightarrow
  r_0 \rightarrow x^{-2} y x y^{-1};\\
14: r_1 \rightarrow  y r_1^{-1} y^{-1} \Rightarrow
  r_1 \rightarrow y x;\\
15: r_0 \rightarrow  r_0 r_1 \Rightarrow
  r_0 \rightarrow x^{-2} y x^2;\\
16: r_0 \rightarrow  x^2 r_0 x^{-2} \Rightarrow
  r_0 \rightarrow y;\\
17: r_1 \rightarrow  r_1^{-1} \Rightarrow
  r_0 \rightarrow x^{-1} y^{-1};\\
18: r_1 \rightarrow  r_1 r_0 \Rightarrow
  r_0 \rightarrow x^{-1};\\
19: r_1 \rightarrow  r_1^{-1} \Rightarrow
  r_0 \rightarrow x;\\
$

To prove that the all other presentations in \ref{hard4} satisfy the
Andrews-Curtis conjecture we show that they are AC-equivalent
to the presentation $G_1 = <x,y;x y x = y x y, x^2 = y^3>$
which was considered above.

For example, let
\[G_3 = <x,y;x^{-1} y^2 x = y^3,x = y x y x^{-1}>\] and
\[r_0 = x^{-1} y^2 x y^{-3}, \; r_1 = x^2 y^{-1} x^{-1} y^{-1}.\]
The following transformations take $G_3$ into the presentation $G_1$: \\
$
1: r_0 \rightarrow  x r_0 x^{-1} \Rightarrow
r_0 \rightarrow  y^2 x y^{-3} x^{-1}; \\
2: r_0 \rightarrow  r_0 r_1 \Rightarrow r_0
\rightarrow  y^2 x y^{-3} x y^{-1} x^{-1} y^{-1};\\
3: r_1 \rightarrow  x y x^{-2} r_1 x^2 y^{-1} x^{-1} \Rightarrow
 r_1 \rightarrow  y^{-1} x^2 y^{-1} x^{-1};\\
4: r_0 \rightarrow  y^{-2} r_0 y^2 \Rightarrow
 r_0 \rightarrow  x y^{-3} x y^{-1} x^{-1} y;\\
5: r_0 \rightarrow  r_0 r_1 \Rightarrow
r_0  \rightarrow  x y^{-3} x y^{-1} x y^{-1} x^{-1};\\
6: r_1 \rightarrow  x y^{-1} x^{-1} r_1 x y x^{-1} \Rightarrow
 r_1 \rightarrow  x y^{-1} x^{-1} y^{-1} x;\\
7:$ Apply automorphism:
\[ x \rightarrow  y^{-1} x^{-1}, y \rightarrow  x y^{-1} x^{-1} \Rightarrow
r_0 \rightarrow  y^2 x^{-3},
r_1 \rightarrow  y x y x^{-1} y^{-1} x^{-1}.
\]

\textit{Remark}. No positive results have been obtained for presentations
of the series (3) with $n >  2$. In all tested cases ($n = 3,4,7,11$)
the least length of a relator was 5 but the total length of the relators never
was reduced.

\section{Particular examples as generic schemes for infinite series}
The purpose of this section is to show that every particular
balanced presentation of the trivial group which satisfies the
Andrews-Curtis conjecture generates an infinite series of
"similar" presentations which are also satisfy the conjecture.

Let $G=<a,b;r(a,b),s(a,b)>$ and $H=<a,b;u(a,b),v(a,b)>$ be
arbitrary presentations. Denote by $G(H)$ the following
presentation
\[G(H)=<a,b;r(u,v),s(u,v)>\]
which is obtained from $G$ by the substitution
$a \rightarrow u$, $b \rightarrow v$ into the relators $r,s$.

\begin{proposition}
Let $G=<a,b;r,s>$ and $H=<a,b;u,v>$ be presentations of the
trivial group. If $G$ satisfies the Andrews-Curtis conjecture,
then the presentation $G(H)$ is AC-equivalent to the presentation
$H$.
\end{proposition}

Indeed, let $t_1,...,t_n$  be a sequence of
transformations $(AC1)-(AC3)$ that bring $G$
to the trivial presentation of the trivial group. Let $s_1,...,s_n$
be a chain of Andrews-Curtis transformations  obtained from $t_1,...,t_n$ as
follows: \\
if $t_i$ is a transformation of the type $(AC1)-(AC2)$ then
$s_i=t_i$;\\
if $t_i$ is of the type $(AC3)$, say $t_i$ replaces $r_j$ by
$w(a,b)r_j w(a,b)^{-1}$, then $s_i$ replaces $r_j$ by
$w(u,v)r_j w(u,v)^{-1}$.\\
Obviously,
\[G(H)s_1...s_n = H.\]

It is clear that if the presentation $H$ satisfies the
Andrews-Curtis conjecture then so does $G(H)$.

As an example, we show that  the following well known trick for producing
presentations of the trivial group which we think is due to B.H.Neumann can
be obtained by the construction above. Suppose that
\[<a,b;r,s>\]
is a presentation of the trivial group. Consider now the
presentation
\[<a,b;r^s=r^2,s^r=s^2>.\]
The above presentation is again that of the trivial group.

The presentation
\[G=<a,b;a^b=a^2,b^a=b^2>\]
is a presentation of the trivial group which satisfies the
Andrews-Curtis conjecture. Indeed, let $r_{0} = b^{-1} a b a^{-2}$ and
$r_{1} = a^{-1} b a b^{-2}$. Apply the following transformation: \\
$
1: r_1 \rightarrow a r_1 a^{-1} \Rightarrow
  r_1 \rightarrow b a b^{-2} a^{-1};\\
2: r_1 \rightarrow  b^{-1} r_1 b \Rightarrow
  r_1 \rightarrow a b^{-2} a^{-1} b;\\
3: r_1 \rightarrow r_1 r_0 \Rightarrow
  r_1 \rightarrow a b^{-1} a^{-2};\\
4: r_0 \rightarrow  a^{-1} r_0 a \Rightarrow
  r_0 \rightarrow a^{-1} b^{-1} a b a^{-1};\\
5: r_0 \rightarrow r_0 r_1 \Rightarrow
  r_0 \rightarrow a^{-1} b^{-1} a^{-1};\\
6: r_0 \rightarrow  a r_0 a^{-1} \Rightarrow
  r_0 \rightarrow b^{-1} a^{-2};\\
7: r_1 \rightarrow  a^{-1} r_1 a \Rightarrow
  r_1 \rightarrow  b^{-1} a^{-1};\\
8: r_1 \rightarrow r_1^{-1} \Rightarrow
  r_1 \rightarrow a b;\\
9: r_1 \rightarrow  r_1 r_0 \Rightarrow
  r_1 \rightarrow a^{-1};\\
10: r_1 \rightarrow  r_1^{-1} \Rightarrow
  r_1 \rightarrow a;\\
11: r_0 \rightarrow  r_0 r_1 \Rightarrow
  r_0 \rightarrow b^{-1} a^{-1};\\
12: r_0 \rightarrow  r_0 r_1 \Rightarrow
  r_0 \rightarrow b^{-1};\\
13: r_0 \rightarrow  r_0^{-1} \Rightarrow
  r_0 \rightarrow b.\\
$

Now let \[H=<a,b;r,s>\] be an arbitrary presentation of the trivial
group. According to the Proposition 1,
\[G(H) = <a,b;r^s=r^2,s^r=s^2>\]
is AC-equivalent to the presentation $H$.

\section{Acknowledgements}

This paper grew out of Gilbert Baumslag's suggestion that genetic algorithms
could be used to good effect in the computational software package MAGNUS.
I would like also to acknowledge the help and numerous suggestions of
Alexei G. Myasnikov and Vladimir Shpilrain during the course of this research.


\end{document}